\documentclass[10pt,a4paper]{article}
\usepackage[mathscr,mathcal]{eucal}
\usepackage{colortbl}
\usepackage{amsfonts}
\usepackage{amssymb}
\usepackage{amsmath}
\usepackage{xypic}
\usepackage{enumerate}
\newtheorem{theorem}{Theorem}[section]
\newtheorem{definition}[theorem]{Definition}
\newtheorem{ypoth}[theorem]{Assumptions}

\newtheorem{proposition}[theorem]{Proposition}

\newtheorem{ex}[theorem]{Example}

\newcommand{\blambda}{\boldsymbol{\lambda}}
\newcommand{\bnu}{\boldsymbol{\nu}}
\newcommand{\bmu}{\boldsymbol{\mu}}

\title{Schur elements and basic sets for cyclotomic Hecke algebras}
\author{Maria Chlouveraki\footnote{University of Edinburgh, e-mail:  maria.chlouveraki@ed.ac.uk\ }\   \ \&   Nicolas Jacon\footnote{Universit\'e Franche-Comt\'e, e-mail: njacon@univ-fcomte.fr}}
\date{}

\begin{document}
\maketitle

\begin{abstract}
We study the Schur elements and the $a$-function for cyclotomic Hecke algebras. 
 As a consequence, we show the existence of canonical basic sets, 
 as defined 
  by Geck-Rouquier,  for certain complex reflection groups. This includes 
  the case of  finite Weyl groups for all choices of parameters (in characteristic $0$).

\end{abstract}

\section{Introduction}
Let $V$ be a finite dimensional complex vector space. A complex reflection group is, by definition,
 a finite subgroup of $\textrm{GL} (V)$ generated by pseudo-reflections (\emph{i.e.,} non-trivial elements 
  acting trivially on a hyperplane). The complex reflection groups
  naturally generalize the finite Weyl groups.  
   Several results and observations suggest that many properties of Weyl groups can be extended to this wider class of   groups.
  In particular,  to any complex reflection group $W$, one can attach a generic  Hecke algebra $\mathcal{H} (W)$ over a ring $A$ which can be seen as a ``quantization" 
  of the group algebra $A[W]$.  
  
  If $K$ is a finite abelian extension of $\mathbb{Q}$ and $q$ is an indeterminate, there is a remarkable specialization of the Hecke algebra  $\mathcal{H} (W)$: the cyclotomic 
specialization $\varphi$. This specialization allows the definition of an algebra which has been intensively studied in the literature:  
the   cyclotomic 
  Hecke algebra
   $\mathcal{H}_{\varphi}$ over $\mathbb{Z}_K[q,q^{-1}]$. 
   One can define a symmetrizing trace form on  $\mathcal{H}_{\varphi}$ which makes $\mathcal{H}_{\varphi}$ into a symmetric algebra.
    As a consequence, to each of the simple   $\mathcal{H}_{\varphi}$-modules we can associate an important polynomial which ``controls"  the trace form and part
     of the representation theory of $\mathcal{H}_{\varphi}$: the Schur element.  The first aim of this paper is to study these elements
      and establish identities between the Schur elements of different cyclotomic Hecke algebras (Propositions  \ref{rho and rho-gamma}, \ref{translation} and \ref{a=A}).

   A natural and important problem in the representation theory of cyclotomic Hecke algebra is to find 
    ``good" parametrizations of the set of simple modules. When the algebra is (split) semisimple, by a 
    well known result (Tits' deformation theorem) 
    this  reduces to the problem of parametrizing the irreducible representations of $W$. The simple modules 
     are then labeled by an indexing set $\Lambda$ of $\textrm{Irr} (W)$.  
   The problem is far more difficult when the algebra is not semisimple.
    Nevertheless, it may be attacked   using the associated decomposition matrix and the theory of basic sets developed by Geck and Rouquier in \cite{GR}. A \emph{basic set} is, by definition, 
    a subset of $\Lambda$ which is in bijection with the set of simple $\mathcal{H}_{\varphi}$-modules.
    This bijection is obtained naturally with the use of the decomposition matrix and an ordering on 
    $\Lambda$. In principle, it is not clear that a basic set always exists. However, 
      the theoretic existence of such a set has been proved when $W$ is a Weyl group
      and the associated weight function is positive  (including the so called ``equal parameter case") and, in certain cases, when $W$ is a complex reflection  group of type $G(r,p,n)$.
     In these cases, the ordering is given by the so called Lusztig's $a$-function,    which can be defined  using the Schur elements.  
   
        Using the study of the Schur elements in the first  part and a careful study of the $a$-function, we will be able to generalize the theorems of existence
         of basic sets to other cases. In particular, we give a complete answer to the question of the existence of basic sets for finite Weyl groups for all choices of parameters
          in characteristic $0$. 
         We then study in detail  the consequences of our results
           in the cases of type $A_{n-1}$, type $B_n$ and then more generally type $G(d,1,n)$ (that is the case of Ariki-Koike algebras).   \\
\\
  {\bf Acknowledgements.} The authors would like to thank Meinolf Geck, Iain Gordon and J\'eremie Guilhot  for useful comments  and discussion.

\section{Schur elements of Hecke algebras}

In this section, we recall the notions of  Schur elements for generic and cyclotomic Hecke algebras and give the definition of the $a$-function. We study what happens to the Schur elements and the $a$-function when we change the parameters that define a cyclotomic Hecke algebra in specific ways. We note that a similar problem has been considered  in the context of Kazhdan-Lusztig cells for finite Weyl groups by C. Bonnaf\'e in \cite{B}.

\subsection{Symmetric algebras}

Let $R$ be an integrally closed Noetherian ring and let $A$ be an $R$-algebra which is free and finitely generated as an $R$-module. We say that $A$ is \emph{symmetric}, if there exists a linear form $t:A \rightarrow R$ which satisfies the following conditions:
\begin{enumerate}[(1)]
\item $t(aa')=t(a'a)$ for all $a,a' \in A$,
\item the bilinear form $\tilde{t}: A \times A \rightarrow R,\,(a,a') \mapsto t(aa')$ is non-degenerate.
\end{enumerate}
The form $t$ is called \emph{symmetrizing} for $A$. 

From now on, let us suppose that the algebra $A$ is symmetric with symmetrizing form $t$. Let $K$ be a field extension of the field of fractions of $R$ such that the algebra $KA:=K \otimes_R A$ is split semisimple.  We denote by $\mathrm{Irr}(KA)$ the set of irreducible characters of the algebra $KA$. Geck has shown (cf.~\cite{Ge}) that $$t=\sum_{\chi \in \textrm{Irr}(KA)}\frac{1}{s_\chi}\chi,$$
where $s_\chi$ is the \emph{Schur element} associated to the
irreducible character $\chi$. The Schur element $s_\chi$ belongs to the integral closure of $R$ in $K$ (cf.~\cite[Proposition 7.3.9]{GePf}) and depends only on the symmetrizing form and the isomorphism class of the representation.

\subsection{Generic Hecke algebras}

Let $K$ be a finite abelian extension of $\mathbb{Q}$ and let $V$ be a finite dimensional
$K$-vector space. Let $W$ be a finite subgroup of $\mathrm{GL}(V)$ generated
by pseudo-reflections and acting irreducibly on $V$. We denote by
$\mathcal{A}$ the set of its reflecting hyperplanes. We set
$V^{\textrm{reg}}:= V-\bigcup_{H \in \mathcal{A}}H$. For $x_0 \in
V^{\textrm{reg}}$, we define $B:=\Pi_1(V^{\textrm{reg}}/W,x_0)$ the
braid group associated with $W$.
For every orbit $\mathcal{C}$ of $W$ on $\mathcal{A}$, let $\textbf{S}_\mathcal{C}$ be the set of the monodromy generators around the
images in $V^{\textrm{reg}}/W$ of the elements of  $\mathcal{C}$. Moreover,
we set
$e_{\mathcal{C}}$ the common order of the subgroups $W_H$, where $H$
is any element of $\mathcal{C}$ and $W_H$ is the pointwise stabilizer of $H$.

We choose a set of indeterminates
$\textbf{u}=(u_{\mathcal{C},j})_{(\mathcal{C} \in
\mathcal{A}/W)(0\leq j \leq e_{\mathcal{C}}-1)}$ and we denote by
$\mathbb{Z}[\textbf{u},\textbf{u}^{-1}]$ the Laurent polynomial ring
in all the indeterminates $\textbf{u}$. We define the \emph{generic Hecke
algebra} $\mathcal{H}(W)$ of $W$ to be the quotient of the group algebra
$\mathbb{Z}[\textbf{u},\textbf{u}^{-1}]B$ by the ideal generated by
the elements of the form
$$(\textbf{s}-u_{\mathcal{C},0})(\textbf{s}-u_{\mathcal{C},1}) \ldots (\textbf{s}-u_{\mathcal{C},e_{\mathcal{C}}-1}),$$
where $\mathcal{C}$ runs over the set $\mathcal{A}/W$ and
$\textbf{s}$ runs over $\textbf{S}_\mathcal{C}$.

We will now make some assumptions for the algebra $\mathcal{H}(W)$ which
have been verified for all but a finite number of irreducible complex reflection groups
(\cite[remarks before 1.17, \S 2]{BMM2},\,\cite{GIM}).

\begin{ypoth}\label{ypo}
The algebra $\mathcal{H}(W)$ is a free
$\mathbb{Z}[\textbf{\emph{u}},\textbf{\emph{u}}^{-1}]$-module of rank equal to the order of $W$. Moreover, there exists a unique symmetrizing form $t$ for $\mathcal{H}(W)$
which has the properties described in \emph{\cite[\S$2$A]{BK}}.
\end{ypoth}

From now on, let us denote by $\mu(K)$ the group of all the roots of unity in $K$
and set $\zeta_{d}:=\mathrm{exp}(2\pi i/d)$ for all $d \geq 1$. 
 Given that the assumptions $\ref{ypo}$ are satisfied, we
 have the following result by Malle (\cite[Theorem 5.2]{Ma4}).

\begin{theorem}\label{Semisimplicity Malle}
Let $\textbf{\emph{v}}=(v_{\mathcal{C},j})_{(\mathcal{C} \in
\mathcal{A}/W)(0\leq j \leq e_{\mathcal{C}}-1)}$ be a set of indeterminates
such that 
$v_{\mathcal{C},j}^{|\mu(K)|}=\zeta_{e_\mathcal{C}}^{-j}u_{\mathcal{C},j}$, for all $\mathcal{C},j$.
Then the $K(\textbf{\emph{v}})$-algebra
$K(\textbf{\emph{v}})\mathcal{H}(W)$ is split semisimple.
\end{theorem}

\subsection{Cyclotomic Hecke algebras}

Let $q$ be an indeterminate and set $y^{|\mu(K)|}:=q$.

\begin{definition}\label{specialization}
A cyclotomic specialization
$\varphi:
\mathbb{Z}_K[\textbf{\emph{v}},\textbf{\emph{v}}^{-1}]\rightarrow
\mathbb{Z}_K[y,y^{-1}]$ is a $\mathbb{Z}_K$-algebra morphism with the following properties:
\begin{enumerate}[(1)]
  \item $\varphi: v_{\mathcal{C},j} \mapsto y^{m_{\mathcal{C},j}}$ where
  $m_{\mathcal{C},j} \in \mathbb{Z}$ for all $\mathcal{C},\,j$,
  \item if $z$ is another
 indeterminate, the element of $\mathbb{Z}_K[y,y^{-1},z]$ defined by\\
  $$\Gamma_\mathcal{C}(y,z):=\prod_{j=0}^{e_\mathcal{C}-1}(z-\zeta_{e_\mathcal{C}}^jy^{m_{\mathcal{C},j}})$$
  is invariant under the action of $\textrm{\emph{Gal}}(K(y)/K(q))$ for all  $\mathcal{C} \in \mathcal{A}/W$.
\end{enumerate}
We usually describe the morphism $\varphi$ by the formula:
$$ u_{\mathcal{C},j} \mapsto \zeta_{e_\mathcal{C}}^jq^{m_{\mathcal{C},j}}.$$ 
\end{definition}

If $\varphi$ is a cyclotomic specialization, then the
corresponding \emph{cyclotomic Hecke algebra} is the
$\mathbb{Z}_K[q,q^{-1}]$-algebra, denoted by $\mathcal{H}_\varphi$,
which is obtained as the specialization of $\mathcal{H}(W)$ via
 $\varphi$. The algebra $\mathcal{H}_\varphi$ has a symmetrizing form $t_\varphi$
defined as  the specialization of the canonical form $t$. We also have an analogue of Theorem
$\ref{Semisimplicity Malle}$ for cyclotomic Hecke algebras (\cite[Proposition 4.3.4]{springer}):

\begin{proposition}\label{semisimplicity chlou}
The $K(y)$-algebra $K(y)\mathcal{H}_\varphi$ is split semisimple.
\end{proposition}

Now, by Tits' deformation theorem (cf., for example, \cite[Theorem 7.4.6]{GePf}), we obtain that the specialization
 $v_{\mathcal{C},j}\mapsto 1$ induces the following bijections between sets of simple modules:
 $$\operatorname{Irr}(K({\bf v}) \mathcal{H}(W)) \leftrightarrow
 \operatorname{Irr}(K(y) \mathcal{H}_{\varphi}) \leftrightarrow
 \operatorname{Irr}(W)$$
 Hence, 
if $\Lambda$ is an indexing set for the irreducible representations of $W$ such that
    $$\operatorname{Irr}(W)=\left\{ V^{\blambda}\ |\ \blambda\in \Lambda\right\},$$
    then we can write:
    $$\operatorname{Irr}(K(y) \mathcal{H}_{\varphi})=\left\{ V_{\varphi}^{\blambda}\ |\ \blambda\in \Lambda\right\}.$$
    \subsection{Functions $a$ and $A$}\label{a}

Following the notations in \cite[6B]{BMM2} for every element $P(y)
\in \mathbb{C}[y,y^{-1}]$, we call
\begin{itemize}
  \item \emph{valuation of $P(y)$ at $y$} and denote by $\mathrm{val}_y(P)$ the order of $P(y)$
  at 0,
  \item \emph{degree of $P(y)$ at $y$} and denote by $\mathrm{deg}_y(P)$ the negative of the
  valuation of $P(1/y)$.
\end{itemize}
Recall that $y^{|\mu(K)|}=q$. We set
$$\mathrm{val}_q(P):=\frac{\mathrm{val}_y(P)}{|\mu(K)|} \,\textrm{ and }\,\,
\mathrm{deg}_q(P):=\frac{\mathrm{deg}_y(P)}{|\mu(K)|}.$$ 

Let $\blambda\in\Lambda$ and denote by $s^{\varphi}_{\blambda}$ the Schur element of $V^{\blambda}_{\varphi}  $ with respect to the symmetrizing form $t_\varphi$. We know that $s^{\varphi}_{\blambda}$ belongs to 
$\mathbb{Z}_K[y,y^{-1}]$. We define
$$a^{\varphi}_{\blambda}:=-\mathrm{val}_q(s_{\blambda}^{\varphi}) \,\textrm{ and }\,
A^{\varphi}_{\blambda}:=\mathrm{deg}_q(s_{\blambda}^{\varphi}).\footnote{In \cite{springer} and \cite{BK}, the $a$-function is defined 
 to be the negative of the one given here. We have chosen here to follow the definition of \cite{GR}.}$$

\subsection{Change of parameters}\label{change}

The generic Hecke algebra $\mathcal{H}(W)$ of $W$ is generated over $\mathbb{Z}[\textbf{u},\textbf{u}^{-1}]$ by a finite number of elements $(\textbf{s}^\mathcal{C}_i)_{1 \leq i \leq n_\mathcal{C}}$, where $\mathcal{C}$ runs over the set $\mathcal{A}/W$, satisfying
\begin{itemize}
\item braid relations,
\item $(\textbf{s}^\mathcal{C}_i-u_{\mathcal{C},0})(\textbf{s}^\mathcal{C}_i-u_{\mathcal{C},1}) \cdots (\textbf{s}^\mathcal{C}_i-u_{\mathcal{C},e_{\mathcal{C}}-1})=0$, for all $\mathcal{C},i$.
\end{itemize}

Let $\textbf{m}=(m_{\mathcal{C},j})_{(\mathcal{C} \in
\mathcal{A}/W)(0\leq j \leq e_{\mathcal{C}}-1)}$ be a family of integers. Consider the cyclotomic specialization $\varphi_{\bf m} :u_{\mathcal{C},j} \mapsto \zeta_{e_\mathcal{C}}^jq^{m_{\mathcal{C},j}}$ and denote by $\mathcal{H}_\textbf{m}$ the corresponding cyclotomic Hecke algebra. Then $\mathcal{H}_\textbf{m}$ is generated by the elements $(\textbf{s}^\mathcal{C}_i)_{(\mathcal{C} \in
\mathcal{A}/W)(1 \leq i \leq n_\mathcal{C})}$, which satisfy the same braid relations as before and
\begin{center}
$(\textbf{s}^\mathcal{C}_i-q^{m_{\mathcal{C},0}})(\textbf{s}^\mathcal{C}_i-\zeta_{e_\mathcal{C}} q^{m_{\mathcal{C},1}}) \cdots (\textbf{s}^\mathcal{C}_i-\zeta_{e_\mathcal{C}} ^{e_{\mathcal{C}}-1}q^{m_{\mathcal{C},e_{\mathcal{C}}-1}})=0$, for all $\mathcal{C},i$.
\end{center}

For each  ${\mathcal{C}} \in \mathcal{A}/W$,  there exists a natural action of the  cyclic group  
$\mathbb{Z}/{e_{\mathcal{C}}}\mathbb{Z}$
 on $\mathbb{Z}^{e_{\mathcal{C}}}$ by cyclic permutation:
 $$\sigma_{\mathcal{C}} : \textbf{m}_\mathcal{C}=\left(m_{\mathcal{C},0},\,m_{\mathcal{C},1},\,\ldots,\,m_{\mathcal{C},e_{\mathcal{C}}-1}\right) \mapsto \sigma_{\mathcal{C}}.\textbf{m}_\mathcal{C}:=
 \left(m_{\mathcal{C},1},\,m_{\mathcal{C},2},\,\ldots,\,m_{\mathcal{C},0}\right).$$
  This induces an action of the group $G:=\prod_{ {\mathcal{C}} \in \mathcal{A}/W} \left(\mathbb{Z}/{e_{\mathcal{C}}}\mathbb{Z}\right)$
   on  the set $\mathcal{M}:=\prod_{ {\mathcal{C}} \in \mathcal{A}/W}  \mathbb{Z}^{e_{\mathcal{C}}}$:
$$\begin{array}{ccc}
G\times \mathcal{M} & \to & \mathcal{M}\\
\left({\bf g}= (g_{\mathcal{C}})_{ {\mathcal{C}} \in \mathcal{A}/W},\, {\bf m}=({\bf m}_{\mathcal{C}})_{ {\mathcal{C}} \in \mathcal{A}/W}\right)& \mapsto & {\bf g}.{\bf m}:=(g_{\mathcal{C}}.{\bf m}_{\mathcal{C}})_{\mathcal{C} \in \mathcal{A}/W}. 
\end{array}$$

Fix  ${\bf g}=(g_{\mathcal{C}})_{ {\mathcal{C}} \in \mathcal{A}/W}\in G$. Abusing notation, we will identify $\mathbb{Z}/{e_{\mathcal{C}}}\mathbb{Z}$ with the set  of its coset representatives $\{0,1,...,e_{\mathcal{C}}-1\}\subset \mathbb{Z}$. 
 For all   ${\mathcal{C}} \in \mathcal{A}/W$ and $i$ such that 
 $1\leq i \leq n_{\mathcal{C}}$,  we set
$$\tilde{\textbf{s}}^{\mathcal{C}}_i:=\zeta_{e_{\mathcal{C}}}^{- g_{\mathcal{C}}  } {\textbf{s}}^{\mathcal{C}}_i.$$ 
The elements of 
$$ \tilde{\textbf{S}}^{\mathcal{C}}:=\bigcup_{\mathcal{C} \in 
\mathcal{A}/W} \bigcup_{  {1 \leq i \leq n_\mathcal{C}} } \left\{\tilde{\textbf{s}}_i^{\mathcal{C}}\right\}$$
 generate $\mathcal{H}_\textbf{m}$. Moreover, for all  ${\mathcal{C}} \in \mathcal{A}/W$ and  $i\,(1\leq i \leq n_{\mathcal{C}})$, the
  generator $\tilde{\textbf{s}}^{\mathcal{C}}_i$ satisfies:
$$(\tilde{\textbf{s}}^{\mathcal{C}}_i-q^{m_{\mathcal{C}, g_{\mathcal{C}}   }})\,(\tilde{\textbf{s}}^{\mathcal{C}}_i-\zeta_{e_{\mathcal{C}}} q^{m_{\mathcal{C},g_{\mathcal{C}} +1}})\,\cdots\,(\tilde{\textbf{s}}^{\mathcal{C}}_i-\zeta_{e_{\mathcal{C}}}^{e_{\mathcal{C}}-1} q^{m_{\mathcal{C},g_{\mathcal{C}}+e_\mathcal{C} -1}})=0,$$
where all subscripts are understood modulo $e_\mathcal{C}$.

Now consider the cyclotomic specialization $\varphi_{{\bf g}. {\bf m}}$ of $\mathcal{H}(W)$. It is  given by:
 $$u_{\mathcal{C},j} \mapsto \zeta_{e_\mathcal{C}}^jq^{m_{\mathcal{C},j+g_{\mathcal{C}}}}
$$ 
 We denote by $\mathcal{H}_{{\bf g}.{\textbf{m}}}$ the corresponding cyclotomic Hecke algebra. 
We can assume that  $\mathcal{H}_{{\bf g}.{\textbf{m}}}$ is generated by the elements of  $\tilde{\textbf{S}}^{\mathcal{C}}$. Then $\mathcal{H}_{{\bf g}.{\textbf{m}}}=\mathcal{H}_{\textbf{m}}$. 
 
 There exists an isomorphism $\Gamma^{\bf g}: \mathcal{H}_{\textbf{m}} \rightarrow \mathcal{H}_{{\bf g}.{\textbf{m}}}$ given by:
 $$\Gamma^{\bf g}(\textbf{s}^{\mathcal{C}}_i)=  \zeta_{e_{\mathcal{C}}}^{g_{\mathcal{C}}} \tilde{\textbf{s}}^{\mathcal{C}}_i$$
 for all  ${\mathcal{C}} \in \mathcal{A}/W$ and $i\,(1\leq i \leq n_{\mathcal{C}})$. 
The isomorphism $\Gamma^{\bf g}$ induces a natural action of
 the group $G$  on the parametrizing set $\Lambda$ as follows: If ${\bf g}\in G$
  and  $\blambda\in \Lambda$, let $\rho^{{\bf g}.{\bf m }}_{\blambda}$ be the associated irreducible representation of $\mathcal{H}_{{\bf g}. {\textbf{m}}}$. 
 The representation  of  $\mathcal{H}_{\textbf{m}}$  defined by the composition $\rho^{{\bf g}. {\bf m }}_{\blambda} \circ \Gamma^{\bf g}$   is  afforded by a simple $\mathcal{H}_{{\textbf{m}}}$-module $V^{\bmu}_{\varphi_{\bf m}}$. We set 
 $$ {\blambda}^{\bf g}:=\bmu$$
 and we obtain a well-defined action: 
 $$\begin{array}{ccc}
 G\times \Lambda & \to & \Lambda\\
 \left({\bf g},\,\blambda \right)& \mapsto & {\blambda}^{\bf g}.
 \end{array}$$
\begin{proposition}\label{rho and rho-gamma}
Let $s^{\varphi_{\bf m}}_{\blambda^{{\bf g}}}$ be the Schur element of
 $\mathcal{H}_{\emph{\textbf{m}}}$ associated to $\rho_{\blambda^{\bf g}}^{\bf m}$ and let $s^{\varphi_{{\bf g}. {\bf m }}}_{\blambda}$ be the Schur element of $\mathcal{H}_{\emph{{\bf g}.\textbf{m}}}$ associated to $\rho^{{\bf g}. {\bf m }}_{\blambda}$. We have
$s^{\varphi_{\bf m}}_{\blambda^{\bf g}}=s^{\varphi_{{\bf g}.{\bf m }}}_{\blambda}$. In particular, $a^{\varphi_{\bf m}}_{\blambda^{\bf g}}=a^{\varphi_{{\bf g}.{\bf m}}}_{\blambda}$.
\end{proposition}
\begin{proof}{Let $\mathcal{H}:=\mathcal{H}_{\textbf{m}}=\mathcal{H}_{{\bf g}.{\textbf{m}}}$. For
 all $\mathcal{C} \in \mathcal{A}/W$, we have
$$\rho^{\bf m }_{\blambda^{\bf g}}\left(\textbf{s}^{\mathcal{C}}_i\right) = \left( \rho^{{\bf g}. {\bf m} }_{\blambda} \circ \Gamma^{\bf g} \right) (\textbf{s}^{\mathcal{C}}_i)=
 \zeta_{e_{\mathcal{C}}} \rho^{{\bf g}. {\bf m}}_{\blambda} (\tilde{\textbf{s}}^{\mathcal{C}}_i)=
 \rho^{{\bf g}. {\bf m}}_{\blambda}  (\textbf{s}^{\mathcal{C}}_i) \,\,\textrm{ for all } i\,(1\leq i \leq n_{\mathcal{C}}).$$
Therefore, the representations  $\rho^{\bf m }_{\blambda^{{\bf g}}}$ and $ \rho^{{\bf g}.{\bf m}}_{\blambda} $ are isomorphic irreducible representations of $\mathcal{H}$. Since the Schur elements depend only on  the isomorphism class of the representation, we deduce that $s^{\varphi_{\bf m}}_{\blambda^{\bf g}}=s^{\varphi_{{\bf g}.{\bf m }}}_{\blambda}$.}
\end{proof}

Given $\blambda \in \Lambda$ and ${\bf g}=(g_{\mathcal{C}})_{ {\mathcal{C}} \in \mathcal{A}/W}\in G$, it is natural to ask how $\blambda^{\bf g}$ is computed\footnote{The authors would like to thank Iain Gordon for pointing that out.}. Since the action of $G$ on $\Lambda$ does not depend on ${\bf m}$, we can assume, without loss of generality,  that
${\bf m}={\bf 0} \in \prod_{ {\mathcal{C}} \in \mathcal{A}/W}  \mathbb{Z}^{e_{\mathcal{C}}}$. Then we have
$\mathcal{H}_{\bf m}=\mathcal{H}_{{\bf g}.{\textbf{m}}} =\mathcal{H}_{\bf 0} \cong \mathbb{Z}_K[W]$ and
$\rho^{\bf 0}_{\blambda^{\bf g}}=\rho^{{\bf 0} }_{\blambda} \circ \Gamma^{\bf g}$. Due to the definition of 
$\Gamma^{\bf g}$, we obtain that 
$$\rho^{\bf 0}_{\blambda^{\bf g}}\,\cong\, \rho^{{\bf 0} }_{\blambda} \otimes \left( \prod_{\mathcal{C} \in \mathcal{A}/W}\left( \mathrm{det}_\mathcal{C} \right)^{g_\mathcal{C}}\right),$$
where $\mathrm{det}_\mathcal{C}$ is the linear character of $W$ defined by
$$\mathrm{det}_\mathcal{C}({\bf s}_i^{\mathcal{C}'})=
\left\{ \begin{array}{ll}
\mathrm{det}_V({\bf s}_i^\mathcal{C})=\zeta_{e_{\mathcal{C}}} &\textrm{ if } \mathcal{C}=\mathcal{C}',\\
&\\
1 & \textrm{ if } \mathcal{C} \neq \mathcal{C}',
\end{array}\right.$$
for all $\mathcal{C},\,\mathcal{C'} \in \mathcal{A}/W$ and $i\,(1 \leq i \leq n_{\mathcal{C}'})$.
Moreover, if $\blambda_{\bf triv} \in \Lambda$ corresponds to the trivial character for $W$, then we have
$$\rho^{\bf 0}_{\blambda_{\bf triv}^{\bf g}}=\prod_{\mathcal{C} \in \mathcal{A}/W}\left( \mathrm{det}_\mathcal{C} \right)^{g_\mathcal{C}}.$$
We deduce that $\blambda^{\bf g}$ is the element of $\Lambda$ which satisfies:
$$\rho^{\bf 0}_{\blambda^{\bf g}}\,\cong\, \rho^{{\bf 0} }_{\blambda} \otimes \rho^{{\bf 0}}_{\blambda_{\bf triv}^{\bf g}}.$$

\begin{ex}\emph{Let $W=G_4$. In this case, there exists a single hyperplane orbit $\mathcal{C}$ under the action of $W$. We write $\chi_{d,b}$ for the irreducible character of $W$ with
dimension $d$ and valuation of fake degree $b$. Hence, the irreducible characters of $G_4$ are the following:
$$\chi_{1,0},\, \,\chi_{1,4}, \,\,\chi_{1,8}, \,\,\chi_{2,5},\,\, \chi_{2,3},\,\, \chi_{2,1},\,\, \chi_{3,2},$$
where $\chi_{1,0}$ is the trivial character for $G_4$ and
$\chi_{1,4}=\mathrm{det}_{\mathcal{C}}$.}

\emph{The generic Hecke algebra $\mathcal{H}(G_4)$ is generated over the Laurent polynomial ring in three indeterminates
$$\mathbb{Z}[u_{\mathcal{C},0},u_{\mathcal{C},0}^{-1},u_{\mathcal{C},1},u_{\mathcal{C},1}^{-1},u_{\mathcal{C},2},u_{\mathcal{C},2}^{-1}]$$
by the elements $\textbf{s}$ and $\textbf{t}$ satisfying the relations:
$$\textbf{sts}=\textbf{tst}\,\,\textrm{ and }\,\,
(\textbf{s}-u_{\mathcal{C},0})(\textbf{s}-u_{\mathcal{C},1})(\textbf{s}-u_{\mathcal{C},2})=
(\textbf{t}-u_{\mathcal{C},0})(\textbf{t}-u_{\mathcal{C},1})(\textbf{t}-u_{\mathcal{C},2})=0.$$
Let $\textbf{m}=(m_{\mathcal{C},0},\,m_{\mathcal{C},1},\,m_{\mathcal{C},2})$ be a family of three integers. We consider the following two cyclotomic specializations of $\mathcal{H}(G_4)$:
$$\varphi_{\bf m}: \left\{ 
\begin{array}{lcl} 
u_{\mathcal{C},0} & \mapsto &q^{m_{\mathcal{C},0}} \\
u_{\mathcal{C},1} & \mapsto &\zeta_3q^{m_{\mathcal{C},1}} \\
u_{\mathcal{C},2} & \mapsto &\zeta_3^2q^{m_{\mathcal{C},2}} \\
\end{array} \right. 
\,\,\textrm{ and }\,\,\,\,
\varphi_{\sigma_\mathcal{C}.{\bf m}}: \left\{ 
\begin{array}{lcl} 
u_{\mathcal{C},0} & \mapsto &q^{m_{\mathcal{C},1}}\\
u_{\mathcal{C},1} & \mapsto &\zeta_3q^{m_{\mathcal{C},2}} \\
u_{\mathcal{C},2} & \mapsto &\zeta_3^2q^{m_{\mathcal{C},0}}. \\
\end{array} \right. 
$$
We denote by $\mathcal{H}_{\textbf{m}}$ and $\mathcal{H}_{\sigma_\mathcal{C}.{\bf m}}$ respectively the corresponding cyclotomic Hecke algebras. If $\chi \in \mathrm{Irr}(W)$, then $s^{\varphi_{\bf m}}(\chi)$ denotes the Schur element of $\chi$ in $\mathcal{H}_{\textbf{m}}$ and $s^{\varphi_{\sigma_\mathcal{C}.{\bf m}}}(\chi)$ denotes the Schur element of $\chi$ in $\mathcal{H}_{\sigma_\mathcal{C}.{\bf m}}$ We obtain the following equalities between the Schur elements of the two algebras:
$$\begin{array}{ccc}
s^{\varphi_\textbf{m}}(\chi_{1,0})=s^{\varphi_{\sigma_\mathcal{C}.{\bf m}}}(\chi_{1,8}),& 
s^{\varphi_\textbf{m}}(\chi_{1,4})=s^{\varphi_{\sigma_\mathcal{C}.{\bf m}}}(\chi_{1,0}),& 
s^{\varphi_\textbf{m}}(\chi_{1,8})=s^{\varphi_{\sigma_\mathcal{C}.{\bf m}}}(\chi_{1,4}), \\  & &\\
 s^{\varphi_\textbf{m}}(\chi_{2,5})=s^{\varphi_{\sigma_\mathcal{C}.{\bf m}}}(\chi_{2,1}),&
 s^{\varphi_\textbf{m}}(\chi_{2,3})=s^{\varphi_{\sigma_\mathcal{C}.{\bf m}}}(\chi_{2,5}), &
 s^{\varphi_\textbf{m}}(\chi_{2,1})=s^{\varphi_{\sigma_\mathcal{C}.{\bf m}}}(\chi_{2,3}), \\  & &\\
 &  s^{\varphi_\textbf{m}}(\chi_{3,2})=s^{\varphi_{\sigma_\mathcal{C}.{\bf m}}}(\chi_{3,2}). & 
 \end{array}$$}
\end{ex}

Now, if $\textbf{n}=(n_{\mathcal{C},j})_{(\mathcal{C} \in \mathcal{A}/W)(0\leq j \leq e_{\mathcal{C}}-1)}
\in \mathcal{M}$, then
$\varphi_{{\bf m}+{\bf n}}$ is the cyclotomic specialization
$$u_{\mathcal{C},j} \mapsto \zeta_{e_\mathcal{C}}^jq^{m_{\mathcal{C},j}+n_{\mathcal{C},j}},$$
such that $\varphi_{\bf m}(v_{\mathcal{C},j})=y^{m_{\mathcal{C},j}}$ and
$\varphi_{{\bf m}+{\bf n}}(v_{\mathcal{C},j})=y^{m_{\mathcal{C},j}+n_{\mathcal{C},j}}$ for all $\mathcal{C},\,j$. 

\begin{proposition}\label{translation}
Let $\blambda \in \Lambda$. 
If ${\bf n} \in \mathcal{M}^G$, i.e., $n_{\mathcal{C},0}=n_{\mathcal{C},1}=\cdots=
n_{\mathcal{C},e_{\mathcal{C}}-1}$ for all $\mathcal{C} \in \mathcal{A}/W$, then
$s_{\blambda}^{\varphi_{\bf m}}=s_{\blambda}^{{\varphi_{{\bf m}+{\bf n}}}}$.
In particular, $a_{\blambda}^{\varphi_{\bf m}}=a_{\blambda}^{{\varphi_{{\bf m}+{\bf n}}}}$.
\end{proposition}
\begin{proof}{Let us denote by $s_{\blambda}^{\textbf{v}}$ the Schur element of $\mathcal{H}(W)$ associated to $\blambda\in \Lambda$. Then $s_{\blambda}^{\varphi_{\bf m}}=\varphi_{\bf m}(s_{\blambda}^{\textbf{v}})$ and $s_{\blambda}^{\varphi_{{\bf m}+{\bf n}}}=\varphi_{{\bf m}+{\bf n}}(s_{\blambda}^{\textbf{v}})$. Following \cite[Theorem 4.2.5]{springer}, $s_{\blambda}^{\textbf{v}}$ is a homogeneous Laurent polynomial of degree $0$ in the indeterminates $(v_{\mathcal{C},j})_{0\leq j \leq e_{\mathcal{C}}-1}$ for all $\mathcal{C} \in \mathcal{A}/W$. Due to the assumption on ${\bf n}$,
it is immediate that 
$\varphi_{\bf m}(s_{\blambda}^{\textbf{v}})=\varphi_{{\bf m}+{\bf n}}(s_{\blambda}^{\textbf{v}})$.
}
\end{proof}

Finally, let $\varphi_{{-\bf m}}$ be the cyclotomic specialization
$u_{\mathcal{C},j} \mapsto \zeta_{e_\mathcal{C}}^jq^{-m_{\mathcal{C},j}}$. The following proposition gives some equalities involving the functions $a$ and $A$.

\begin{proposition}\label{a=A}
Let $\blambda \in \Lambda$. We have $a_{\blambda}^{\varphi_{\bf m}}=A_{\blambda}^{{\varphi_{-{\bf m}}}}$ and $\,A_{\blambda}^{\varphi_{\bf m}}=a_{\blambda}^{{\varphi_{-{\bf m}}}}$.
\end{proposition}
\begin{proof}{We have $s_{\blambda}^{\varphi_{\bf m}}=\varphi_{\bf m}(s_{\blambda}^{\textbf{v}})$ and $s_{\blambda}^{\varphi_{-{\bf m}}}=\varphi_{-{\bf m}}(s_{\blambda}^{\textbf{v}})$, where $s_{\blambda}^{\textbf{v}}$ denotes the Schur element of $\mathcal{H}(W)$ associated to $\blambda\in \Lambda$. Therefore, $s_{\blambda}^{\varphi_{-{\bf m}}}(y)=s_{\blambda}^{\varphi_{{\bf m}}}(y^{-1})$.
The result follows immediately from the definition of the functions $a$ and $A$, given in \S \ref{a}.}
\end{proof}

\section{Basic Sets}

We now study some of the consequences of Proposition \ref{rho and rho-gamma}. We focus in particular on the modular representation theory of Hecke algebras. 
 From known results and the above ones, we will be able to generalize the notion of ``canonical basic sets" as defined by Geck and Rouquier in \cite{GR}.

\subsection{Decomposition maps}\label{decomp}

Let $\theta:\mathbb{Z}_K[y,y^{-1}]\to k$ be a ring homomorphism 
such that $k$ is the field of fractions of $\theta (\mathbb{Z}_K[y,y^{-1}])$ and
$\theta (q )=\xi$, where
    $\xi\in k^{\times}$ is an element which has
a  root of order $|\mu (K)|$ in $k^{\times}$. Considering $k$ as a $\mathbb{Z}_K[y,y^{-1}]$-module via $\theta$, we assume that the algebra
$k\mathcal{H}_{\varphi}:=k\otimes_{\mathbb{Z}_K[y,y^{-1}] }\mathcal{H}_{\varphi}$ is split.  

As noted above, we 
 have a canonical way to parametrize the simple modules for $K (y)\mathcal{H}_{\varphi}$. It is also desirable 
  to obtain a ``good" parametrization of the simple $k\mathcal{H}_{\varphi}$-modules. As $k\mathcal{H}_{\varphi}$ is 
   not semisimple in general, Tits' deformation theorem cannot be applied. However, following 
   \cite[\S 7.4]{GePf}, one can use the associated decomposition
    matrix to solve that problem. 
    
    Let 
$R(K(y)\mathcal{H}_{\varphi})$ (respectively $R(k\mathcal{H}_{\varphi})$) be the Grothendieck group
 of finitely generated $K(y)\mathcal{H}_{\varphi}$-modules (respectively $k\mathcal{H}_{\varphi}$-modules). 
  It is generated by the classes $[U]$ of the simple  $K(y)\mathcal{H}_{\varphi}$-modules (respectively $k\mathcal{H}_{\varphi}$-modules) $U$.
   Then we obtain a well-defined  decomposition map
$$d^{\varphi}_{\theta}:R_0 (K(y)\mathcal{H}_{\varphi})    \to R_0 (k\mathcal{H}_{\varphi})$$
such that, for all $\lambda\in \Lambda$ , we have
$$d^{\varphi}_{\theta}([V^{\blambda}_{\varphi}])=\sum_{M\in \operatorname{Irr}(k\mathcal{H}_{\varphi})}
 [ V_{\varphi}^{\blambda} : M] [M].$$
The matrix
\[ D^{\varphi}_{\theta}=\left([ V_{\varphi}^{\blambda} :M]\right)_{\blambda\in \Lambda,\,M \in \operatorname{Irr}(k \mathcal{H}_{\varphi} ) }\]
is called the \emph{decomposition matrix associated with} $\theta$. 
    
    We now take the setting of \S \ref{change}. Let ${\bf g} \in G=\prod_{ {\mathcal{C}} \in \mathcal{A}/W} \left( \mathbb{Z}/{e_{\mathcal{C}}}\mathbb{Z} \right)$
   and ${\bf m} \in \mathcal{M}=\prod_{ {\mathcal{C}} \in \mathcal{A}/W}  \mathbb{Z}^{e_{\mathcal{C}}}$. Then 
    we have two cyclotomic specializations  $\varphi_{\bf m}$ and $\varphi_{{\bf g}. \textbf{m} }$
     (with corresponding cyclotomic Hecke algebras $\mathcal{H}_{\bf m}$ and $\mathcal{H}_{{\bf g}. \textbf{m}}$ respectively).
 The isomorphism $\Gamma^{\bf g}:\mathcal{H}_{\bf m}\to \mathcal{H}_{{\bf g}. {\bf m}}$ induces two group isomorphisms at the level of the Grothendieck group of finitely 
  generated modules:
   $$[\Gamma^{\bf g}]_K: R_0 (K(y)\mathcal{H}_{{\bf g}. {\bf m}}) \to R_0(K(y)\mathcal{H}_{{\bf m}} )$$
 and 
   $$[\Gamma^{\bf g}]_k: R_0 (k\mathcal{H}_{{\bf g}. {\bf m}} )\to R_0(k\mathcal{H}_{{\bf m}}).$$
As a consequence, we obtain the following commutative diagram:\\

 \centerline{
\xymatrix{
    R_0 (K(y)\mathcal{H}_{{\bf g}. {\bf m}})      \ar[r]^{[\Gamma^{\bf g}]_K} \ar[d]^{d^{\varphi_{{\bf g}. {\bf m}} }_{\theta}} & R_0 (K(y)\mathcal{H}_{\bf m})    \ar[d]^{d^{\varphi_{\bf m}}_{\theta}} \\
   R_0 (k\mathcal{H}_{{\bf g}. {\bf m}} )    \ar[r]_{[\Gamma^{\bf g}]_k} & R_0 (k\mathcal{H}_{\bf m}) . 
  }
}$ $\\
\noindent Thus, for all $\blambda\in \Lambda$ and $M\in \operatorname{Irr}( k\mathcal{H}_{{{\bf g}. {\bf m}}})$, we have 
$$[V_{\varphi_{{\bf g}. {\bf m}}}^{\blambda}:M] =[V_{\varphi_{\bf m}}^{\blambda^{{\bf g}}}:M^{\Gamma_k^{\bf g}}] ,$$
where $M^{\Gamma^{\bf g}_k}$ denotes the simple $k\mathcal{H}_{{\bf m}}$-module obtained by composition with $\Gamma^{\bf g}_k : k\mathcal{H}_{{\bf m}}\to k\mathcal{H}_{{{\bf g}. {\bf m}} }$

The notion of basic sets of simple modules for Hecke algebras has first been considered by Geck in \cite{cells} when $W$ is a finite Weyl group.
 We recall this definition of these objects  in the following section.

\subsection{Parametrizing Hecke algebra modules}
 The theory of basic sets gives canonical ways to parametrize the simple modules for Hecke algebras. 
  In this section, we assume that we have a cyclotomic Hecke algebra $\mathcal{H}_{\varphi}$ and a 
specialization    $\theta:\mathbb{Z}_K[y,y^{-1}]\to k$ as in \S \ref{decomp}. By the previous section, we have a decomposition matrix:
$$ D^{\varphi}_{\theta}=\left([ V_{\varphi}^{\blambda} :M]\right)_{\blambda\in \Lambda,\,M \in \operatorname{Irr}(k \mathcal{H}_{\varphi} ) }.$$
\begin{definition}\label{basique}
We say that $\mathcal{H}_{\varphi}$  admits a basic set $\mathcal{B}(\varphi)\subset \Lambda$ with respect to
$\theta:\mathbb{Z}_K[y,y^{-1}]\to k$    if and only if the following two conditions are satisfied:
  \begin{enumerate}[(1)]
\item   For all $M\in\operatorname{Irr}(k \mathcal{H}_{\varphi})$, there exists $\blambda_M\in \mathcal{B}(\varphi)$ such that
 $$[V_{\varphi}^{\blambda_M}:M] =1\textrm{ and }a^{\varphi}_{\bmu}>a^{\varphi}_{\blambda_M}\textrm{ if }[V_{\varphi}^{\bmu}:M] \neq 0.$$
\item The map 
 $$\begin{array}{cll}
 \operatorname{Irr}(k \mathcal{H}_{\varphi})   &\to& \mathcal{B}(\varphi)\\
M &\mapsto& \lambda_M
\end{array}$$
is a bijection. \end{enumerate}
\end{definition}

Assume that $k \mathcal{H}_{\varphi}$   admits a basic set $\mathcal{B}(\varphi)\subset \Lambda$ with respect to
  $\theta$. This implies that the associated decomposition matrix has a lower 
  triangular shape with one's along the diagonal for a ``good" order on  $\Lambda$ induced by the map $a^{\varphi}$. Hence, 
   it gives a way to label $\operatorname{Irr}(k \mathcal{H}_{\varphi})$.

  The existence of basic sets has been shown for all Hecke algebras of finite Weyl groups in the case where 
   the associated ``weight function" is positive (\cite{Gsurvey}), for special cases of Ariki-Koike algebras (\cite{Jkyoto}),
    and more generally for certain cyclotomic Hecke algebras of type $G(r,p,n)$ (\cite{JG}). 
   We will generalize these results in the following section. Before this, we need to study the consequences of a change of parameters for the basic sets.

\begin{proposition}\label{existence}
We keep the notation of  \S\ref{change} and assume that $k \mathcal{H}_{\bf m}$ admits a basic set $\mathcal{B}(\varphi_{\bf m})\subset \Lambda$ with respect to $\theta:\mathbb{Z}_K[y,y^{-1}]\to k$.
Then,  for all ${\bf g}\in G$, the algebra $k \mathcal{H}_{{\bf g}. {\bf m}}$ 
 admits a basic set $\mathcal{B}(\varphi_{{\bf g}. {\bf m }})\subset \Lambda$ with respect to $\theta$. Moreover, 
 $$\mathcal{B}(\varphi_{{\bf g}.{\bf m}})=\left\{ \blambda^{{\bf g}^{-1}}\ |\ \blambda\in \mathcal{B}(\varphi_{\bf m})\right\}.$$ 
\end{proposition}
\begin{proof}{ Let $M\in \textrm{Irr}(k\mathcal{H}_{{\bf g}.{\bf m}})$ and set $N:= M^{\Gamma^{\bf g}_k}\in   \textrm{Irr}(k\mathcal{H}_{\bf m})$.  
By assumption, there exists 
 $\blambda_{N} \in\mathcal{B}(\varphi_{\bf m})$ such that $ [V_{\varphi_{\bf m}}^{\blambda_{N}}:N] =1$
  and 
$ [V_{\varphi_{\bf m}}^{\bmu}:N] \neq 0$ implies 
$a_{\bmu}^{\varphi_{\bf m}}>a_{\blambda_{N}}^{\varphi_{\bf m}}$.

Now, combining Proposition \ref{rho and rho-gamma} with the results of the  previous section, we obtain that
 $\bnu_M:=\blambda_{N}^{{\bf g}^{-1}}$ satisfies  
   $ [V_{\varphi_{{\bf g}.{\bf m}}}^{\bnu_{M}}:M] =1$
  and that
$ [V_{\varphi_{{\bf g}.{\bf m}}}^{\bmu}:M] \neq 0$ implies 
$a_{\bmu}^{\varphi_{{\bf g}.{\bf m}}}>a_{\bnu_{M}}^{\varphi_{{\bf g}.{\bf m}}}$. As the map
$$\begin{array}{ccc}
\Lambda&\to &\Lambda\\
 \blambda&\mapsto &\blambda^{\bf {g}^{-1}}\end{array}$$
  is a bijection, the proposition follows.}
\end{proof}

Similarly to the above proposition, we have the following result 
 which follows from Proposition \ref{translation}. 

\begin{proposition}\label{remarksuperbedemariacequinel'empechepasdetrepenible}
 Assume that $k \mathcal{H}_{\bf m}$ admits a basic set $\mathcal{B}(\varphi_{\bf m})\subset \Lambda$ with respect to $\theta:\mathbb{Z}_K[y,y^{-1}]\to k$.
 Then, for all  ${\bf n} \in \mathcal{M}^G$, the algebra  
  $k \mathcal{H}_{{\bf m}+{\bf n}}$ 
 admits a basic set $\mathcal{B}(\varphi_{{\bf m }+{\bf n}})\subset \Lambda$ with respect to $\theta$. Moreover,
 $$\mathcal{B}(\varphi_{{\bf m}+{\bf n}})= \mathcal{B}(\varphi_{\bf m}).$$ 

\end{proposition}

In the following sections, we study consequences of the above results on the basic sets 
 for Hecke algebras of finite Weyl groups and, more generally, complex reflection groups.

\section{Applications}

Using the above results, we will be able to generalize the theorems of existence of basic sets for Hecke algebras of finite Weyl groups. 

\subsection{Basic sets for Hecke algebras of finite Weyl groups}
 Hecke algebras of finite Weyl groups  are special cases of the Hecke algebras of complex reflection groups. 
 In these cases,  for all  $\mathcal{C}\in \mathcal{A}/W$, we have $e_{\mathcal{C}}=2$ and thus, for all $i\,(1\leq i \leq n_{\mathcal{C}})$, the generator $\textbf{s}_i^{\mathcal{C}}$ satisfies
$$(\textbf{s}_i^{\mathcal{C}}-q^{m_{\mathcal{C},0}})(\textbf{s}_i^{\mathcal{C}}+q^{ {m_{\mathcal{C},1}}})=0,$$
 for some integers $m_{\mathcal{C},0}$,  $m_{\mathcal{C},1}$.
 
 The following theorem shows the existence of a canonical basic set for all finite Weyl groups and all choices of ${\bf m}$ in characteristic $0$.
 It  is an extension of a result   originally obtained by M.~Geck and R.~Rouquier (\cite{GR}) 
  in the special case  where there exists $a\in \mathbb{N}$ such that 
 $m_{\mathcal{C},0}=a$  and  $m_{\mathcal{C},1}=0$ for all $\mathcal{C}\in \mathcal{A}/W$.
 Some generalization of the result has  been considered by M.~Geck (\cite{Gsurvey}) and M.~Geck and the second author (\cite{GJ}) in the case where  $m_{\mathcal{C},0}\in \mathbb{N}$  and  $m_{\mathcal{C},1}=0$, 
   for all $\mathcal{C}\in \mathcal{A}/W$. 
  We refer to \cite{GJlivre} for a complete survey of this theory.  
  The proof of the following generalization uses as main ingredients these 
  results together with the properties of the $a$-function established in the previous sections.
 \begin{theorem}\label{main}
Assume that $W$ is a finite Weyl group and that the characteristic of $k$ is zero. Then
 for all ${\bf m}=(m_{\mathcal{C},j})_{\mathcal{C}\in \mathcal{A}/W,\, j=0,1}$, the algebra 
 $\mathcal{H}_{\bf m}$ admits a canonical basic set in the sense of Definition \ref{basique} with respect to any specialization. 
 \end{theorem}
 \begin{proof}{
If  $m_{\mathcal{C},0}\in \mathbb{N}$  and  $m_{\mathcal{C},1}=0$ for all $\mathcal{C}\in \mathcal{A}/W$, then, as noted above,   the result has already  been proved. 
 Using Proposition \ref{remarksuperbedemariacequinel'empechepasdetrepenible}, we can conclude in the case where we have  
  $m_{\mathcal{C},0}\geq m_{\mathcal{C},1}$ for all  $\mathcal{C}\in \mathcal{A}/W$. Finally, the general case follows from
  Proposition \ref{existence}.
   }
 \end{proof}
 
 Using the results in \cite{Gsurvey} and the same proof as above, we have an analogue of the above theorem in arbitrary (good) characteristic, assuming that 
 the conjectures  ${\bf P1-P15}$ of \cite{Lu} hold.
  In the next section, we study in detail the case of the Ariki-Koike algebras and in particular, the 
   question of the parametrization of the basic sets in type $A_{n-1}$.

\subsection{Ariki-Koike algebras}

The group $G(d,1,n)$ is the group of all monomial $n \times n$ matrices with entries in $\mathbb{Z}/d\mathbb{Z}$. It is isomorphic to the wreath product $(\mathbb{Z}/d\mathbb{Z}) \wr
\mathfrak{S}_n$  and its field of definition is $K:=\mathbb{Q}(\zeta_d)$. 
The \emph{cyclotomic Ariki-Koike algebra} of $G(d,1,n)$ (cf. \cite{ArKo}, \cite{BM}) is the algebra  generated over $\mathbb{Z}_K[q,q^{-1}]$ 
by the elements $\mathrm{\textbf{s}},\mathrm{\textbf{t}}_1,\mathrm{\textbf{t}}_2,\ldots,\mathrm{\textbf{t}}_{n-1}$ satisfying the relations:
\begin{itemize}
\item $\mathrm{\textbf{s}}\mathrm{\textbf{t}}_1\mathrm{\textbf{s}}\mathrm{\textbf{t}}_1=\mathrm{\textbf{t}}_1\mathrm{\textbf{s}}\mathrm{\textbf{t}}_1\mathrm{\textbf{s}}$, $\mathrm{\textbf{s}}\mathrm{\textbf{t}}_j=\mathrm{\textbf{t}}_j\mathrm{\textbf{s}} \textrm{ for } j\neq 1$,
\item $\mathrm{\textbf{t}}_j\mathrm{\textbf{t}}_{j+1}\mathrm{\textbf{t}}_j=\mathrm{\textbf{t}}_{j+1}\mathrm{\textbf{t}}_j\mathrm{\textbf{t}}_{j+1}$,  $ \mathrm{\textbf{t}}_i\mathrm{\textbf{t}}_j=\mathrm{\textbf{t}}_j\mathrm{\textbf{t}}_i \textrm{ for } |i-j|>1$,
\item $(\mathrm{\textbf{s}}-q^{m_{\mathcal{C},0}})(\mathrm{\textbf{s}}-\zeta_d q^{m_{\mathcal{C},1}})\ldots(\mathrm{\textbf{s}}-\zeta_d^{d-1}q^{m_{\mathcal{C},d-1}})=(\mathrm{\textbf{t}}_j-q^{m_{\mathcal{C'},0}}   )(\mathrm{\textbf{t}}_j+q^{m_{\mathcal{C'},1}})=0$.
\end{itemize}

It is well-known that, in this case, the simple modules are naturally parametrized by the set of 
 $d$-partitions of $n$: 
$$\Pi^d (n):=\left\{ (\lambda^{(0)},\,\ldots,\,\lambda^{(d-1)}) \,\left|\right.\,\forall i\in \{0,\,\ldots,\,d-1\}\,\, \lambda^{(i)}\in \Pi ({n_i}),\,\, \sum_{k=0}^{d-1}{n_k}=n\right\},$$
where $\Pi (m)$ denotes the set of partitions of rank $m$ (see \cite{ArKo}).   
The question of the existence of basic sets for Ariki-Koike algebras has been studied in \cite{Jkyoto} and \cite{Juglov} (see also \cite{GJ} for the case $W=G(r,p,n)$). 
In these papers, the explicit parametrization of the basic sets by subsets of $\Pi^d (n)$ has been established.

Using Proposition \ref{existence}, we obtain the existence and the explicit parametrization of the basic sets in other cases. 
 For this,  we give here the action of the group $G=\mathbb{Z}/d\mathbb{Z}\times \mathbb{Z}/2\mathbb{Z}$ on $\Pi^d (n)$ by describing the action of its generators. We also 
  give the identities on the $a$-function that Proposition \ref{rho and rho-gamma} implies. 
  We keep the notation of \S \ref{change}.
  
\begin{itemize}
\item  If ${\bf g}=(1,0)\in G$, then 
  $$(\lambda^{(0)},\,\lambda^{(1)},\,\ldots,\,\lambda^{(d-1)})^{{\bf g}}=(\lambda^{(d-1)},\,\lambda^{(0)},\,\ldots,\,\lambda^{(d-2)})$$
  and we have
  $$a^{\varphi_{\bf m}}_{(\lambda^{(0)},\,\lambda^{(1)},\,\ldots,\lambda^{(d-1)})}=a^{\varphi_{{\bf g}. {\bf m }}}_{(\lambda^{(1)},\,\ldots,\,\lambda^{(d-1)},\,\lambda^{(0)})}.$$
  \item If ${\bf g}=(0,1)\in G$, then 
  $$(\lambda^{(0)},\,\lambda^{(1)},\,\ldots,\,\lambda^{(d-1)})^{{\bf g}}=({\lambda^{(0)}} ',\,{\lambda^{(1)}} ',\,\ldots,\,{\lambda^{(d-1)}}'),$$
  where $'$ denotes the usual conjugation of partitions. Moreover, we have
    $$a^{\varphi_{\bf m}}_{(\lambda^{(0)},\,\lambda^{(1)},\,\ldots,\,\lambda^{(d-1)})}=a^{\varphi_{{\bf g}. {\bf m }}}_{({\lambda^{(0)}} ',\,{\lambda^{(1)}} ',\,\ldots,\,{\lambda^{(d-1)}}')}.$$
\end{itemize}

  Here we only wish to study in detail the case of finite Weyl groups   where the existence of canonical basic sets has been established for all choices 
  of parameters  by Theorem \ref{main} (which is not the case when $d>2$).
We focus in particular on type $A_{n-1}$.

\subsection{An example : type $A_{n-1}$} 
Let us now study the case of type $A_{n-1}$.  The associated cyclotomic Hecke algebra is  generated over $\mathbb{Z}[q,q^{-1}]$ 
by the elements $\mathrm{\textbf{t}}_1,\mathrm{\textbf{t}}_2,\ldots,\mathrm{\textbf{t}}_{n-1}$ satisfying the relations:
\begin{itemize}
\item $\mathrm{\textbf{t}}_j\mathrm{\textbf{t}}_{j+1}\mathrm{\textbf{t}}_j=\mathrm{\textbf{t}}_{j+1}\mathrm{\textbf{t}}_j\mathrm{\textbf{t}}_{j+1}$,  $ \mathrm{\textbf{t}}_i\mathrm{\textbf{t}}_j=\mathrm{\textbf{t}}_j\mathrm{\textbf{t}}_i \textrm{ for } |i-j|>1$,
\item $(\mathrm{\textbf{t}}_j-q^{m_{\mathcal{C},0}}   )(\mathrm{\textbf{t}}_j+q^{m_{\mathcal{C},1}})=0$.
\end{itemize}
The set  $\Lambda$ is the set $\Pi (n)$ of partitions of rank $n$ and $G$ is the group $\mathbb{Z}/2\mathbb{Z}$. 
  The bijection on $\Lambda$ induced by the non-trivial element of $G$ is then the conjugation of partitions.  
   Assume that $m_{\mathcal{C},0}\in \mathbb{N}$
 and $m_{\mathcal{C},1}=0$. Then, if $\lambda=(\lambda_1,\,\ldots,\,\lambda_r) \in \Pi (n)$, we have 
 $$a^{\varphi_{\bf m}}_{\lambda}=m_{\mathcal{C},0} \sum_{i=1}^r (i-1)\lambda_i.$$ 
Following the discussion in the previous section, we deduce that, in the general case,
$$a^{\varphi_{\bf m}}_{\lambda}=\left\{ \begin{array}{ll}
 (m_{\mathcal{C},0}-m_{\mathcal{C},1}) \displaystyle{\sum_{i=1}^r (i-1)\lambda_i}& \textrm{ if }m_{\mathcal{C},0}\geq m_{\mathcal{C},1},\\
 (m_{\mathcal{C},1}-m_{\mathcal{C},0}) \displaystyle{\sum_{i=1}^s (i-1)\lambda_i'}& \textrm{ if }m_{\mathcal{C},0} < m_{\mathcal{C},1},
 \end{array}\right.$$
where $\lambda'=(\lambda_1',\,\ldots,\,\lambda_s ')$ denotes the conjugate partition of $\lambda$ (hence, $s=\lambda_1$). 
Proposition \ref{a=A} implies that $A^{\varphi_{\bf m}}_\lambda=a^{\varphi_{\bf m}}_{\lambda'}$.
 Obviously, if $m_{\mathcal{C},0}\geq m_{\mathcal{C},1}$, we have 
 $$\lambda \unrhd \mu \Rightarrow a^{\varphi_{\bf m}}_\lambda\leq a^{\varphi_{\bf m}}_\mu,$$
 where $\unrhd$ denotes the usual dominance order on partitions. 
   The next result follows from \cite{DJ0} and Proposition \ref{existence}.
\begin{proposition}
Let $\theta$ be a specialization 
 such that $\theta (q)^{m_{\mathcal{C},0}-m_{\mathcal{C},1}}$ is a primitive root of $1$ of order $e>1$.
\begin{enumerate}[(1)]
\item Assume that ${\bf m}$ is such that $m_{\mathcal{C},0}> m_{\mathcal{C},1}$.   
Then $\mathcal{H}_{\bf m}$ admits a basic set $\mathcal{B}(\varphi_{\bf m})\subset \Lambda$ with respect to
  $\theta$  and we have
  $$\mathcal{B}(\varphi_{\bf m})=\operatorname{Reg}_e (n)$$
  where the set $\operatorname{Reg}_e (n)$ of $e$-regular partitions is defined by 
  $$\lambda\notin \operatorname{Reg}_e (n)\iff \exists i\in \mathbb{N},\ \lambda_i=\cdots=\lambda_{i+e-1}\neq 0.$$
\item Assume that ${\bf m}$ is such that $m_{\mathcal{C},0}< m_{\mathcal{C},1}$.  Then $\mathcal{H}_{\bf m}$ admits a basic set $\mathcal{B}(\varphi_{\bf m})\subset \Lambda$ with respect to
  $\theta$  and we have
  $$\mathcal{B}(\varphi_{\bf m})=\operatorname{Res}_e (n)$$
  where the set $\operatorname{Res}_e (n)$ of $e$-restricted partitions is defined by 
  $$\lambda\in \operatorname{Res}_e (n)\iff \forall i\in \mathbb{N},\ \lambda_i-\lambda_{i+1}\leq e-1.$$
  \end{enumerate}
\end{proposition}

In type $B_n$, the parametrization of the basic sets is more complicated to describe and 
 uses objects coming from the crystal basis theory for the quantum group $\mathcal{U}_q (\widehat{\mathfrak{sl}}_e)$. A complete study of this case, including connections with other 
  remarkable objects related to Hecke algebras (the constructible representations), will be presented in \cite{Jcons}.

\end{document}